\documentclass[12pt,a4paper,british,english]{amsart}
\usepackage[T1]{fontenc}
\usepackage[latin9]{inputenc}
\usepackage{color}
\usepackage{babel}

\usepackage{array}
\usepackage{units}
\usepackage{amsthm}
\usepackage{amsbsy}
\usepackage{amstext}
\usepackage{amssymb}
\usepackage[unicode=true, pdfusetitle,
 bookmarks=true,bookmarksnumbered=false,bookmarksopen=false,
 breaklinks=false,pdfborder={0 0 0},backref=false,colorlinks=true]
 {hyperref}

\makeatletter

\pdfpageheight\paperheight
\pdfpagewidth\paperwidth

\providecommand{\tabularnewline}{\\}

\numberwithin{equation}{section}
\numberwithin{figure}{section}

\usepackage{mathrsfs}
\usepackage{mathtools}
\usepackage{euscript}
\usepackage{upgreek}

\makeatother

\begin{document}
\global\long\def\dual#1{#1^{{\scriptscriptstyle \vee}}}

\global\long\def\trans#1{^{t}\!#1}

\global\long\def\set#1#2{\left\{  #1\, |\, #2\right\}  }

\global\long\def\map#1#2#3{#1\!:#2\!\rightarrow\!#3}

\global\long\def\aut#1{\mathrm{Aut}\!\left(#1\right)}

\global\long\def\End#1{\mathrm{End}\!\left(#1\right)}

\global\long\def\id#1{\mathbf{id}_{#1}}

\global\long\def\ip#1{\left[#1\right]}

\global\long\def\uhp{\mathbf{H}}

\global\long\def\sm#1#2#3#4{\left(\begin{smallmatrix}#1  &  #2\cr\cr#3  &  #4\end{smallmatrix}\right)}

\global\long\def\cyc#1{\mathbb{Q}\left[\zeta_{#1}\right]}

\global\long\def\ZN#1{\left(\mathbb{Z}/#1\mathbb{Z}\right)^{\times}}

\global\long\def\Mod#1#2#3{#1\equiv#2\, \left(\mathrm{mod}\, \, #3\right)}

\global\long\def\psl{\mathrm{PSL}_{2}\!\left(\mathbb{Z}\right)}

\global\long\def\SL{\mathrm{SL}_{2}\!\left(\mathbb{Z}\right)}

\global\long\def\gl#1#2{\mathrm{GL}_{#1}\!\left(#2\right)}

\global\long\def\tr#1{\mathrm{Tr}#1}

\global\long\def\mr#1#2#3#4{\rho\!\left(\!\begin{array}{cc}
 #1  &  #2\\
#3  &  #4\end{array} \!\right)}

\global\long\def\FA{\mathbb{X}}

\global\long\def\FB{\mathbf{\boldsymbol{\Lambda}}}

\global\long\def\FC{\mathsf{h}_{{\scriptscriptstyle 0}}}

\global\long\def\FD{\boldsymbol{\Gamma}}

\global\long\def\FE{\EuScript{P}_{\FB}}

\global\long\def\FF{\mathbb{C}\!\left[J\right]}

\global\long\def\FG{\boldsymbol{\Gamma}}

\global\long\def\sig#1{\upsigma\!\left(#1\right)}

\global\long\def\wof#1#2{\mathcal{M}_{#1}\!\left(#2\right)}

\global\long\def\sof{\mathcal{M}_{0}}

\global\long\def\mof#1#2{\mathsf{M}_{#1}\!\left(#2\right)}

\global\long\def\muf#1#2{\mathsf{M}_{#1}^{\circ}\!\left(#2\right)}

\global\long\def\cof#1#2{\mathsf{S}_{#1}\!\left(#2\right)}

\global\long\def\cuf#1#2{\mathsf{S}_{#1}^{\circ}\!\left(#2\right)}

\global\long\def\kan{\varkappa}

\global\long\def\un{\varsigma(\tau)}

\global\long\def\ws{\boldsymbol{\varpi}}

\global\long\def\dll#1{\varGamma_{#1}}

\global\long\def\ks{\boldsymbol{\uplambda}_{{\scriptscriptstyle +}}}

\global\long\def\kso{\boldsymbol{\uplambda}_{{\scriptscriptstyle -}}}

\global\long\def\gam#1{\upgamma_{#1}}

\global\long\def\mcr#1{\boldsymbol{X}\!\!\left(#1\right)}

\global\long\def\Gmof#1#2{\mathfrak{M}_{#1}\!\left(#2\right)}
\global\long\def\Gcof#1#2{\mathfrak{S}_{#1}\!\left(#2\right)}

\title{The dimension of vector-valued modular forms of integer weight}

\author{P. Bantay}

\curraddr{Institute for Theoretical Physics, Eötvös Loránd University, Budapest}

\email{bantay@poe.elte.hu}

\thanks{Work supported by grant OTKA 78005.}
\begin{abstract}
We present a dimension formula for spaces of vector-valued modular
forms of integer weight in case the associated multiplier system has
finite image, and discuss the weight distribution of the module generators
of holomorphic and cusp forms, as well as the duality relation between
cusp forms and holomorphic forms for the contragredient.
\end{abstract}
\maketitle

\section{Introduction}

The classical theory of scalar modular forms \cite{apostol2,Knopp1970,koblitz}
has been a major theme of mathematics in the last centuries. Its applications
are numerous, ranging from number theory to topology and mathematical
physics, a showpiece being the mathematics involved in the proof of
Fermat's Last Theorem \cite{diamond}. An important tool in the applications
of the theory is the explicit description of the different spaces
of modular forms, which allows to identify highly transcendental functions
via their analytic and transformation properties. In particular, a
major result is the dimension formula for spaces of holomorphic and
cusp forms, which allows to determine explicit bases for these spaces,
and describe the involved algebraic structures very precisely \cite{Lang,Serre}.

While the need for a theory of vector-valued forms, i.e. holomorphic
maps from the complex upper half-plane into a linear space that transform
according to some nontrivial (projective) representation of the modular
group $\SL$, has been recognized long ago, its systematic development
has begun only recently \cite{Knopp2004,Bantay2006,Mason2007,Bantay2008}.
The importance of vector-valued modular forms for mathematics lies,
besides the intrinsic interest of the subject, in the fact that important
classical problems may be reduced to the study of suitable vector-valued
forms, like the theory of Jacobi forms \cite{Eichler1985} or of scalar
modular forms for finite index subgroups \cite{Selberg1965}; from
a modern perspective, trace functions of vertex operator algebras
\cite{FLM1,Kac} satisfying suitable restrictions also provide important
examples of vector-valued modular forms \cite{Zhu1996}. From the
point of view of theoretical physics, vector-valued modular forms
play an important role in string theory \cite{GSW,Polch} and two-dimensional
conformal field theory \cite{DiFrancesco-Mathieu-Senechal}, as the
basic ingredients (chiral blocks) of torus partition functions and
other correlators. 

The above connections justify amply the interest in obtaining explicit
expressions for the dimension of spaces of vector-valued modular forms.
Such results do indeed exist in the literature \cite{Eholzer-Skoruppa,Manschot},
mostly based%
\footnote{Except for \cite{Bantay2009}, which anticipates the results of the
present paper.%
} on the pioneering work \cite{Skoruppa1984}; while the latter, relying
on the Eichler-Selberg trace formula, provides a closed expression
for the dimensions, it doesn't give a constructive procedure for determining
explicit bases, which is a serious drawback from the point of view
of many applications. The present paper offers an alternative approach,
based on the results of \cite{Bantay2006,Bantay2008}, for computing
the dimension of various spaces of vector-valued modular forms, which
is conceptually simpler, and can be modified easily to provide effective
procedures for computing explicit bases.

\section{Vector-valued modular forms\label{sec:general}}

Let $V$ denote a finite dimensional linear space, $\rho\!:\!\FD\!\rightarrow\!\gl{}V$
a representation of $\FD\!=\!\SL$ on $V$, and $w$ an integer. A
(vector-valued) modular form of weight $w$ with multiplier $\rho$
is a map $\FA\!:\!\uhp\!\rightarrow\! V$ that is holomorphic everywhere
in the upper half-plane $\uhp\!=\!\set{\tau}{\mathrm{Im}\tau\!>\!0}$,
and transforms according to the rule\begin{equation}
\FA\!\left(\frac{a\tau+b}{c\tau+d}\right)=\left(c\tau+d\right)^{w}\mr abcd\FA\!\left(\tau\right)\:\label{eq:modtrans}\end{equation}
for $\sm abcd\!\in\!\FD$. 

A form is called weakly holomorphic if it has at worst finite order
poles in the limit $\tau\!\rightarrow\!\mathsf{i}\infty$, i.e. its
Puisseux-expansion in terms of the local uniformizing parameter $q\!=\!\exp\!\left(2\pi\mathsf{i}\tau\right)$
has only finitely many terms with negative exponents; it is holomorphic,
respectively a cusp form if it is bounded (resp. vanishes) as $\tau\!\rightarrow\!\mathsf{i}\infty$,
meaning that its Puisseux-expansion contains only non-negative (resp.
positive) powers of $q$. We'll denote by $\wof w{\rho}$ the (in
general infinite dimensional) linear space of weakly holomorphic forms,
and by $\mof w{\rho}$ and $\cof w{\rho}$ the subspaces of holomorphic
and cusp forms; clearly, we have the inclusions $\cof w{\rho}\!<\!\mof w{\rho}\!<\!\wof w{\rho}$.
An obvious but important observation is that \begin{equation}
\wof w{\rho_{1}\!\oplus\!\rho_{2}}=\wof w{\rho_{1}}\oplus\wof w{\rho_{2}}\:\label{eq:decomp}\end{equation}
for any two representations $\rho_{1}$ and $\rho_{2}$, and a similar
decomposition holds for the spaces of holomorphic and cusp forms,
which implies that \begin{equation}
\begin{aligned}\dim\mof k{\oplus_{i}\rho_{i}} & ={\textstyle \sum}_{i}\dim\mof k{\rho_{i}}\,,\\
\dim\cof k{\oplus_{i}\rho_{i}}\,\, & ={\textstyle \sum}_{i}\dim\cof k{\rho_{i}}\,.\end{aligned}
\label{eq:dimforred}\end{equation}

Note that we recover the classical theory of (scalar) modular forms
of $\SL$ when $\rho\!=\!\rho_{0}$ is the trivial (identity) representation.
In this case the weight should be an even integer for nontrivial forms
to exist, which is non-negative (resp. positive) for holomorphic (resp.
cusp forms). By a well known result \cite{apostol2,Serre}, the spaces
$\mof{2k}{\rho_{0}}$ and $\cof{2k}{\rho_{0}}$ are all finite dimensional:
$\mof 0{\rho_{0}}$ consists of constants, $\mof 2{\rho_{0}}$ is
empty, while $\mof 4{\rho_{0}}$ and $\mof 6{\rho_{0}}$, each having
dimension $1$, are spanned by the Eisenstein series\begin{align}
E_{4}\!\left(q\right)=\, & 1+240\sum\limits _{n=1}^{\infty}\sigma_{3}\left(n\right)q^{n}\,\label{eq:eisexpl4}\\
\intertext{and}E_{6}\!\left(q\right)=\, & 1-504\sum\limits _{n=1}^{\infty}\sigma_{5}\left(n\right)q^{n}\,,\label{eq:eisexpl6}\end{align}
where $\sigma_{k}\!\left(n\right)=\sum_{d|n}d^{k}\,$ is the $k^{\mathrm{th}}$
power sum of the divisors of $n$. What is more, any holomorphic form
may be expressed uniquely as a bivariate polynomial in the Eisenstein
series $E_{4}\!\left(q\right)$ and $E_{6}\!\left(q\right)$, in other
words \begin{equation}
\mathsf{M}=\bigoplus_{k=0}^{\infty}\mof k{\rho_{0}}=\mathbb{C}\!\left[E_{4},E_{6}\right]\,\label{eq:modring}\end{equation}
as graded rings. On the other hand, there are no cusp forms of weight
less than $12$, while $\cof{12}{\rho_{0}}$ is spanned by the discriminant
form \begin{equation}
\Delta\!\left(q\right)\,=\!\frac{1}{1728}\left(E_{4}\!\left(q\right)^{3}-E_{6}\!\left(q\right)^{2}\right)\!=\! q\prod_{n=1}^{\infty}\left(1-q^{n}\right)^{24}\,,\label{eq:deltadef}\end{equation}
and any cusp form of weight $k\!\geq\!12$ is the product of $\Delta\!\left(q\right)$
with a holomorphic form of weight $k\!-\!12$. Finally, the ring $\wof 0{\rho_{0}}$
of scalar weakly-holomorphic forms of weight $0$, which we shall
denote simply by $\sof$ in the sequel, coincides with the univariate
polynomial algebra $\FF$ generated by the Hauptmodul (the trace function
of the Moonshine module \cite{Borcherds1,FLM1}) \begin{equation}
J\!\left(q\right)\!=\!\dfrac{E_{4}\!\left(q\right)^{3}}{\Delta\!\left(q\right)}-744\!=\! q^{-1}\!+\!196884q\!+\cdots\,,\label{eq:Jdef}\end{equation}
and each $\wof{2k}{\rho_{0}}$ is a module over $\sof$ generated
by a single element.

Since multiplying a holomorphic form $\FA\!\left(\tau\right)\!\in\!\mof k{\rho}$
with a scalar holomorphic form $f\!\left(\tau\right)\!\in\!\mof{2n}{\rho_{0}}$
results in a new holomorphic form $f\!\left(\tau\right)\!\FA\!\left(\tau\right)\!\in\!\mof{k+2n}{\rho}$,
and the same is true for cusp forms, the direct sums $\mof{}{\rho}\!=\!\oplus_{k=0}^{\infty}\mof k{\rho}$
and $\cof{}{\rho}\!=\!\oplus_{k=0}^{\infty}\cof k{\rho}$ are (graded)
modules over the ring $\mathsf{M}$ of holomorphic scalar modular
forms. An important result \cite{Marks2009} states that these are
free modules of rank $d$. An interesting question in this respect
is to determine the weight distribution of a set of free generators,
which may be answered by considering the Hilbert-Poincaré series of
these modules \cite{Eisenbud}, i.e. the generating functions $\Gmof{\rho}z\!=\!\sum_{k}\dim\mof k{\rho}z^{k}$
and $\Gcof{\rho}z\!=\!\sum_{k}\dim\cof k{\rho}z^{k}$: the number
of independent generators of weight $k$ equals the coefficient of
$z^{k}$ in the finite polynomials $\left(1\!-\! z^{4}\right)\!\left(1\!-\! z^{6}\right)\Gmof{\rho}z$
and $\left(1\!-\! z^{4}\right)\!\left(1\!-\! z^{6}\right)\Gcof{\rho}z$.

We'll call a representation $\rho\!:\!\FD\!\rightarrow\!\gl{}V$ even
in case $\rho\!\sm{\textrm{-}1}00{\textrm{-}1}\!=\!\id V$, and odd
if $\rho\!\sm{\textrm{-}1}00{\textrm{-}1}\!=\!-\id V$. Any representation
$\rho$ may be decomposed uniquely into a direct sum $\rho\!=\!\rho_{+}\!\oplus\!\rho_{-}$
of even and odd representations, and any irreducible representation
is either even or odd. It follows from Eq.\eqref{eq:modtrans} that
for an even (resp. odd) representation $\rho$ there are no nontrivial
forms of odd (resp. even) weight. Combining this result with Eq.\eqref{eq:decomp},
one gets at once that \begin{equation}
\dim\mof k{\rho}\!=\!\begin{cases}
\dim\mof k{\rho_{+}} & \textrm{if }k\textrm{ is even,}\\
\dim\mof k{\rho_{-}} & \textrm{if }k\textrm{ is odd,}\end{cases}\label{eq:oddeven}\end{equation}
and a similar result for cusp forms. This result shows that it is
enough to treat separately purely even and odd representations, the
general case can be reduced to these.

Since the discriminant form $\Delta\!\left(\tau\right)$ does not
vanish on the upper half-plane \cite{apostol2}, its 12th root $\un\!=\! q{}^{\nicefrac{1}{12}}\prod_{n=1}\left(1\!-\! q^{n}\right)^{2}$
(the square of Dedekind's eta function) is well-defined and holomorphic
on $\uhp$, with an algebraic branch point at the cusp $\tau\!=\!\mathsf{i}\infty$.
Moreover, $\un$ is a weight 1 cusp form with multiplier $\kan$,
where $\kan$ denotes the one dimensional representation of $\SL$
for which%
\footnote{$\kan$ generates the group of linear characters of $\SL$, which
is cyclic of order 12; moreover, $\kan$ is an odd representation,
and tensoring with $\kan$ takes an even representation into an odd
one and \emph{vice versa.} %
}

\begin{equation}
\begin{aligned}\kan\!\sm 0{\textrm{-}1}1{\,0} & \!=-\mathsf{i}\\
\kan\!\sm{\,0}{\,\textrm{-}1}{\,1}{\,\textrm{-}1} & \!=\exp\!\left(\!\dfrac{4\pi\mathsf{i}}{3}\right)\:.\end{aligned}
\label{eq:canonical}\end{equation}
It does follow that, for any representation $\rho$ and any form $\FA\!\in\!\wof w{\rho}$,
one has $\un^{k}\FA\!\left(\tau\right)\!\in\!\wof{w+k}{\rho\!\otimes\!\kan^{k}}$
for all integers\emph{ }$k\!\in\!\mathbb{Z}$; in other words, one
has an injective map\begin{align}
\ws_{k}:\wof w{\rho}\rightarrow & \wof{w+k}{\rho\!\otimes\!\kan^{k}}\label{eq:wsdef}\\
\FA\!\left(\tau\right)\mapsto & \un^{k}\FA\!\left(\tau\right)\,.\nonumber \end{align}
The map $\ws_{k}$ relates forms of different weights with a slightly
different multiplier. Note that, since $\un\!\in\!\cof 1{\kan}$ is
a cusp form, multiplication by a positive power of $\un$ takes a
holomorphic form into a cusp form, i.e. $\ws_{k}\!\left(\mof w{\rho}\right)\!<\!\cof{w+k}{\rho\!\otimes\!\kan^{k}}$
for $k\!>\!0$, in particular \begin{equation}
\dim\mof w{\rho}\!\leq\!\dim\cof{w+k}{\rho\!\otimes\!\kan^{k}}\,.\label{eq:dimbound}\end{equation}

The idea underlying most of what follows is that the injectivity of
the weight-shifting map $\ws_{k}$ allows to reduce the study of forms
of arbitrary weights to that of forms of weight $0$. For example,
the space $\mof w{\rho}$ of weight $w$ holomorphic forms may be
characterized through its image $\muf w{\rho\!\otimes\!\kan^{\textrm{-}w}}\!=\!\ws_{\textrm{-}w}\!\left(\mof w{\rho}\right)\!<\!\wof 0{\rho\!\otimes\!\kan^{\textrm{-}w}}$.
This means that the basic objects of study are the spaces $\wof 0{\rho}$,
together with their various subspaces. Fortunately, the structure
of $\wof 0{\rho}$ is pretty well understood and under control: let's
shortly review the relevant results.

To start with, we make an important restriction on the representation
$\rho$: from now on, we require that $\rho$ has finite image (equivalently,
that its kernel be of finite index). While this requirement could
seem too restrictive at first sight, it is satisfied in most cases
of interest: just to cite an example relevant to physics, the representation
describing the modular properties of the chiral characters of a Rational
Conformal Field Theory has a kernel of finite index \cite{Bantay2003b}.
Moreover, while one may develop a theory for more general representations
$\rho$, both the formulation and the proof of the relevant results
becomes much more cumbersome, and lacks the elegance of the results
to be presented.

The benefits of requiring the image of $\rho$ to be finite are numerous,
the most important being:
\begin{enumerate}
\item $\rho$ is completely reducible (by Maschke's theorem), i.e. it can
be decomposed into a direct sum of irreducible representations;
\item the operator $\rho\!\sm 1101$ can be diagonalized, and its eigenvalues
are roots of unity (since it is of finite order);
\item the kernel of $\rho$ uniformizes a finite sheeted cover of the modular
curve $\uhp/\FD\cong\mathbb{CP}^{1}$.
\end{enumerate}
Let's now turn to the properties of $\wof 0{\rho}$ (recall that $\rho$
is supposed to be even and of finite image). The basic observation
is that the product $f\!\left(\tau\right)\FA\!\left(\tau\right)$
of a weakly holomorphic form $\FA\!\left(\tau\right)\!\in\!\wof 0{\rho}$
with a scalar form $f\!\left(\tau\right)\!\in\!\sof$ is again a weakly
holomorphic form belonging to $\wof 0{\rho}$: in other words, $\wof 0{\rho}$
is an $\sof$-module, what is more, it is a torsion free module. Taking
into account the fact that $\sof$ is the univariate polynomial algebra
$\FF$ generated by the Hauptmodul, this means that actually $\wof 0{\rho}$
is a free module \cite{Eisenbud}, whose rank equals the dimension
$d$ of the representation $\rho$. This means that there exists forms
$\FA_{1},\ldots,\FA_{d}\!\in\!\wof 0{\rho}$ that freely generate
$\wof 0{\rho}$ as an $\sof$-module, i.e. any weakly holomorphic
form $\FA\!\in\!\wof 0{\rho}$ may be decomposed uniquely into a sum
\begin{equation}
\FA\!\left(\tau\right)\!=\!\sum_{i=1}^{d}\wp_{i}\!\left(\tau\right)\FA_{i}\!\left(\tau\right)\,,\label{eq:xdecomp}\end{equation}
where the coefficients $\wp_{1},\ldots,\wp_{d}\!\in\!\sof$ are weight
$0$ weakly holomorphic scalar forms: since $\sof\!=\!\FF$, the coefficients
may be considered as univariate polynomials in the Hauptmodul $J\!\left(\tau\right)$. 

Since, by assumption, $\rho\!\sm 1101$ can be diagonalized, there
exists a diagonalizable operator $\FB$ (called the exponent matrix),
such that \begin{equation}
\rho\!\sm 1101\!=\!\exp\!\left(2\pi\mathsf{i}\FB\right)\,.\label{eq:lambdadef}\end{equation}
Note that $\FB$ is far from unique, its eigenvalues being only determined
up to integers. Taking into account the transformation rule Eq.\eqref{eq:modtrans},
it is clear that for all $\FA\!\in\!\wof 0{\rho}$ the expression
$\exp\!\left(\textrm{-}2\pi\mathsf{i}\FB\tau\right)\FA\!\left(\tau\right)$
is periodic in $\tau$ with period $1$, hence%
\footnote{Recall that $q\!=\!\exp\!\left(2\pi\mathsf{i}\tau\right)$ is the
uniformizing parameter at $\tau\!=\!\mathsf{i}\infty$.%
} \begin{equation}
q^{\textrm{-}\FB}\FA\!\left(q\right)\!=\!\sum_{n}\FA\!\left[n\right]q^{n}\label{eq:qexpan}\end{equation}
for some coefficients $\FA\!\left[n\right]\!\in\! V$, with only finitely
many negative powers of $q$ on the right hand side of Eq.\eqref{eq:qexpan}.
The sum \begin{equation}
\FE\FA\!=\!\sum_{n<0}\FA\!\left[n\right]q^{n}\label{eq:ppartdef}\end{equation}
of these negative powers is the ($\FB$-)principal part of the form
$\FA\!\left(\tau\right)$; clearly, it depends on the actual choice
of $\FB$. 

An important result is that one may always choose $\FB$ such that
the corresponding principal part map $\FE$ is bijective \cite{Bantay2008},
i.e. any form is uniquely determined by its principal part, and any
sum $\sum_{n<0}\FA\!\left[n\right]q^{n}$ is the principal part of
some form $\FA\!\in\!\wof 0{\rho}$. A necessary condition for the
bijectivity of $\FE$ is the relation \begin{equation}
\tr{\,\FB}\!=\! d-\dfrac{\alpha}{2}-\dfrac{\beta_{1}+2\beta_{2}}{3}\,\,,\label{eq:traceformula}\end{equation}
where the integers $\alpha$, $\beta_{1}$ and $\beta_{2}$ are important
numerical characteristics of the representation $\rho$, termed collectively
its signature: $\alpha$ denotes the multiplicity of $-1$ as an eigenvalue
of $\rho\!\sm 0{\textrm{-}1}10$ (which is an involution, since $\rho$
is even), while $\beta_{1}$ and $\beta_{2}$ denote the multiplicities
of $\exp\!\left(\frac{2\pi\mathsf{i}}{3}\right)$ and $\exp\!\left(\frac{4\pi\mathsf{i}}{3}\right)$
as eigenvalues of $\rho\!\sm 0{\textrm{-}1}1{\textrm{-}1}$. Note
that the signature can be determined through the relations\begin{equation}
\begin{aligned}\tr{\,\rho\!\sm 0{\textrm{-}1}10} & =\, d-2\alpha\:,\\
\tr{\,\rho\!\sm 0{\textrm{-}1}1{\textrm{-}1}} & =\, d-\frac{3}{2}\left(\beta_{1}+\beta_{2}\right)+\mathsf{i}\frac{\sqrt{3}}{2}\left(\beta_{1}-\beta_{2}\right)\:.\end{aligned}
\label{eq:sigtrace}\end{equation}

\section{The general dimension formula}

Let's consider an even irreducible representation $\rho\!:\!\SL\!\rightarrow\!\gl{}V$
having finite image, an integer $k$ and an exponent matrix $\FB$
for which the principal part map $\FE$ is bijective; in particular,
$\FB$ has to satisfy the trace formula Eq.\eqref{eq:traceformula}.
Let $\ip x$ denote the integer part of $x\!\in\!\mathbb{R}$, i.e.
the largest integer not exceeding $x$, and for a (diagonalizable)
operator $A$, let $\tr{\ip A}$ denote the sum of the integer parts
of its eigenvalues. The basic dimension formula, from which all others
follow, reads\begin{equation}
\begin{aligned}\dim\mof k{\rho\!\otimes\!\kan^{k}} & =\max\!\left(0,\tr{\ip{\FB\!+\!\dfrac{k}{12}}}\,\right)\\
\dim\,\cof{k\,}{\rho\!\otimes\!\kan^{k}} & =\max\!\left(0,-\tr{\ip{1\!-\!\FB\!-\!\dfrac{k}{12}}}\,\right)\:.\end{aligned}
\label{eq:dimfor1}\end{equation}

Let's see how the above result comes about. The first observation
is that, thanks to the injectivity of the weight shifting map $\ws_{-k}$,
one has $\dim\mof k{\rho\!\otimes\!\kan^{k}}\!=\!\dim\muf k{\rho}$
and $\dim\cof k{\rho\!\otimes\!\kan^{k}}\!=\!\dim\cuf k{\rho}$. By
definition, $\muf k{\rho}$ (resp. $\cuf k{\rho}$) consists of those
weakly holomorphic forms $\FA\!\left(\tau\right)\!\in\!\wof 0{\rho}$
for which $\un^{k}\FA\!\left(\tau\right)$ remains bounded (resp.
vanishes) as $\tau\!\rightarrow\!\mathsf{i}\infty$, i.e. for which
\begin{equation}
q^{\frac{k}{12}+\FB+n}\FA\!\left[n\right]\label{eq:holo1}\end{equation}
tends to a finite limit (resp. vanishes) for each $n$ as $q\!\rightarrow\!0$,
cf. Eq.\eqref{eq:qexpan}. Since the principal part map $\FE$ is
bijective by assumption, any form $\FA\!\in\!\wof 0{\rho}$ is completely
determined by its expansion coefficients $\FA\!\left[n\right]$ with
$n\!<\!0$; in particular, for $n\!\geq\!0$ the expansion coefficients
$\FA\!\left[n\right]$ are linear expressions in the coefficients
$\FA\!\left[m\right]$ with $m\!<\!0$. Choosing a basis in which
$\FB$ is diagonal (this is always possible, thanks to our assumptions
on $\rho$), and denoting by $\FA_{i}$ the component of $\FA$ corresponding
to the eigenvalue $\FB_{i}$, the condition for $\FA$ belonging to
$\muf k{\rho}$ (resp. $\cuf k{\rho}$) is that $\FA_{i}\!\left[n\right]\!=\!0$
\foreignlanguage{british}{provided $\frac{k}{12}+\FB_{i}+n$ is negative
(resp. non-positive).} Let's observe that these conditions constitute
a linear system of equations in the variables $\FA_{j}\!\left[m\right]$
with $m\!<\!0$, and holomorphic (resp. cusp) forms are in one-to-one
correspondence with solutions of this system. \global\long\def\kl{\upmu_{i}}

Let's consider the quantity $\kl\!=\!\frac{k}{12}+\FB_{i}$. If $\kl\!>\!0$,
there are exactly $\left[\kl\right]$ (resp. $-\!\ip{1\!-\!\kl}$)
negative integers $n$ for which $\kl+n$ is non-negative (resp. positive),
and for these values of $n$ the corresponding $\FA_{i}\!\left[n\right]$
may be nonvanishing according to the above. On the other hand, for
$\kl\!\leq\!0$ not only the $\FA_{i}\!\left[n\right]$-s with $n\!<\!0$,
but also the first $-\!\left[\kl\right]$ (resp. $\ip{1\!-\!\kl}$)
components with $n\!\geq\!0$ have to vanish: the later, being linear
expressions in the $\FA_{j}\!\left[m\right]$-s with $m\!<\!0$, supply
us with $-\!\left[\kl\right]$ (resp. $\ip{1\!-\!\kl}$) linear relations
on the coefficients of a holomorphic (resp. cusp) form. Subtracting
the total number $r$ of relations from the number $m$ of possible
nonvanishing coefficients, we get a total of $m\!-\! r\!=\!\sum_{i}\left[\frac{k}{12}\!+\!\FB_{i}\right]\!=\!\tr{\ip{\FB\!+\!\tfrac{k}{12}}}$
(resp. $m\!-\! r\!=\!-\!\tr{\ip{1\!-\!\FB\!-\!\tfrac{k}{12}}}$) free
coefficients by the above reasoning. If $m\!<\! r$, i.e. there are
more relations than nonvanishing coefficients, then the resulting
linear system is overdetermined, having no nontrivial solutions at
all, and $\dim\muf k{\rho}\!=\!0$ (resp. $\dim\cuf k{\rho}\!=\!0$).
On the other hand, for $m\!\geq\! r$ there are less relations than
variables, and $\dim\muf k{\rho}\!=\! m\!-\! r$ (resp. $\dim\cuf k{\rho}\!=\! m\!-\! r$).
Putting all this together, we get Eq.\eqref{eq:dimfor1}. Note that
the restriction to irreducible $\rho$ is important, since the above
argument assumes that the representation space $V$ is a minimal (nontrivial)
invariant subspace for $\rho$; it could very well happen for a reducible
representation that $m\!\geq\! r$, but for some subrepresentation
there are more relations than nonvanishing coefficients, with the
result that the relevant dimension is strictly less than $m\!-\! r$.

While Eq.\eqref{eq:dimfor1} solves the original problem, it still
needs some elaboration. Indeed, one would like a formula expressing
$\dim\mof k{\rho}$ and $\dim\cof k{\rho}$ as a function of $k$
for fixed $\rho$. Of course, it is trivial to arrive to such an expression
from Eq.\eqref{eq:dimfor1}, by simply replacing the representation
$\rho$ with $\rho\!\otimes\!\kan^{{\scriptscriptstyle -}k}$ (note
that, since we have defined exponent matrices for even representations
only, this makes sense for even $k$ only in case $\rho$ is even,
and for odd $k$ only if $\rho$ is odd), giving\begin{equation}
\begin{aligned}\dim\mof k{\rho} & =\max\!\left(0,\tr{\ip{\FB^{\left(k\right)}\!+\!\dfrac{k}{12}}}\right)\,,\\
\dim\,\cof k{\rho}\, & =\max\!\left(0,-\tr{\ip{1\!-\!\FB^{\left(k\right)}\!-\!\dfrac{k}{12}}}\,\right)\:,\end{aligned}
\label{eq:dimfor2}\end{equation}
with $\FB^{\left(k\right)}$ denoting the exponent matrix of $\rho\!\otimes\!\kan^{{\scriptscriptstyle -}k}$. 

The problem with Eq.\eqref{eq:dimfor2} is twofold. First, it makes
reference to the exponent matrix $\FB^{\left(k\right)}$ of the representation
$\rho\!\otimes\!\kan^{{\scriptscriptstyle -}k}$, but one would like
to dispense of the need to compute these quantities, and express the
relevant traces solely in terms of some simple numerical characteristics
of $\rho$. The second problem with Eq.\eqref{eq:dimfor2} is that
it is only valid for irreducible representations, and one would like
a general result valid for any representation. While there is an obvious
solution to this, exploiting the additivity Eq.\eqref{eq:dimforred}
of dimensions and the fact that representations with finite image
are completely reducible, this approach requires the knowledge of
the irreducible decomposition of $\rho$, while one would like an
explicit expression for the dimensions in terms of some global characteristics
of the representation. Since the way to achieve the above goals differs
slightly for even and odd representations, we shall treat these cases
separately in the subsequent sections, starting with the even case.

\section{Even representations\label{sec:even}}

Let's fix an even representation $\rho$, and recall that in this
case there are no forms of odd weight. Our starting point is the observation
that, by the very definition of $\FB^{\left(k\right)}$, \begin{align}
\exp(2\pi\mathsf{i}\FB^{\left(2k\right)})= & \left(\rho\!\otimes\!\kan^{{\scriptscriptstyle -}2k}\right)\!\!\sm 1101=\label{eq:lambdak}\\
\qquad=\kan\!\sm 1101^{{\scriptscriptstyle -}2k}\!\rho\!\sm 1101= & \exp\!\left(2\pi\mathsf{i}\left(\FB\!-\!\tfrac{k}{6}\right)\!\right)\,,\nonumber \end{align}
from which one concludes that all the eigenvalues of 

\begin{equation}
\dll k\!=\!\FB^{\left(2k\right)}\!-\!\FB\!+\!\dfrac{k}{6}\label{eq:Gammakdef}\end{equation}
are necessarily integers, consequently\begin{equation}
\begin{aligned}\tr{\!\ip{\FB^{\left(2k\right)}\!+\!\frac{k}{6}}} & \!=\!\tr{\ip{\FB\!+\!\dll k}}\!=\!\tr{\!\ip{\FB}}\!+\!\tr{\ip{\dll k}}\,,\\
\tr{\!\ip{1\!-\!\FB^{\left(2k\right)}\!-\!\dfrac{k}{6}}} & \!=\!\tr{\ip{1\!-\!\FB\!-\!\dll k}}\!=\!\tr{\!}\ip{1\!-\!\FB}\!-\!\tr{\ip{\dll k}}\,.\end{aligned}
\label{eq:dimfor3}\end{equation}

Upon introducing (for arbitrary even $\rho$) the notations \begin{equation}
\:\quad\begin{aligned}\ks\!\left(\rho\right) & =\!\tr{\!\ip{\FB}}\,,\\
\kso\!\left(\rho\right) & =\!-\!\tr{\!\ip{1\!-\!\FB}}\,\end{aligned}
\label{eq:lamdef}\end{equation}
and \begin{equation}
\gam k\!\left(\rho\right)\!=\!\tr{\!\ip{\dll k}},\label{eq:gamdef}\end{equation}
the formula for even irreducible $\rho$ takes on the form\begin{equation}
\begin{aligned}\dim\mof{2k}{\rho} & =\max\!\left(0,\ks\!+\!\gam k\,\right)\,,\\
\dim\cof{2k}{\rho}\, & =\max\!\left(0,\kso\!+\!\gam k\,\right)\,.\end{aligned}
\label{eq:dimfor4}\end{equation}

A few comments are in order at this point. First, let's note that
the difference $\ks-\kso$ equals the number of integer eigenvalues%
\footnote{This follows from the observation that, for any $x\!\in\!\mathbb{R}$,
the sum $\ip x+\ip{1\!-\! x}$ equals $1$ if $x$ is an integer,
and $0$ otherwise.%
} of $\FB$, i.e. the number of invariant vectors of the operator $\rho\!\sm 1101$;
in particular, one has $\kso\!\leq\!\ks$, in complete accord with
the inclusion $\cof{2k}{\rho}\!<\!\mof{2k}{\rho}$. Moreover, the
integer sequence of $\gam k$-s follows a simple repetitive pattern,
namely (recall that $d$ denotes the dimension of $\rho$) \begin{equation}
\gam{k+6}\!=\!\gam k\!+\! d\,\label{eq:gamperiod}\end{equation}
for all $k$, as a consequence of the fact that the 12th power of
$\kan$ is the identity representation $\rho_{0}$, hence $\FB^{\left(k+12\right)}\!=\!\FB^{\left(k\right)}$.
This means that, since $\gam 0\!=\!0$, the whole sequence is determined
by $\gam 1,\ldots,\gam 5$, and one has $\dim\mof{2k+12}{\rho}\!=\!\dim\mof{2k}{\rho}+d$
provided $\dim\mof{2k}{\rho}\!>\!0$, with a similar result for cusp
forms. Finally, since $\dll k\!\left(\rho\!\otimes\!\kan^{{\scriptscriptstyle -}2n}\right)\!=\!\FB^{\left(k+n\right)}\!-\!\FB^{\left(n\right)}\!+\!\tfrac{k}{6}\!=\!\dll{k+n}\!\left(\rho\right)\!-\!\dll n\!\left(\rho\right)$,
one concludes that \begin{equation}
\gam k\!\left(\rho\!\otimes\!\kan^{{\scriptscriptstyle -}2n}\right)\!=\!\gam{k+n}\!\left(\rho\right)\!-\!\gam n\!\left(\rho\right)\,.\label{eq:gamdif}\end{equation}

For Eq.\eqref{eq:dimfor4} to be effective, it remains to give a practical
method to determine the $\gam k$-s. This is based on the observation
that, because the eigenvalues of $\dll k$ are integers, one has \begin{equation}
\tr{\!\ip{\dll k}}\!=\!\tr{\!\left(\dll k\right)}\!=\!\tr{\FB^{\left(2k\right)}}-\tr{\FB}+\dfrac{kd}{6}\,.\label{eq:trgam1}\end{equation}

But the traces appearing in this expression may be expressed, thanks
to the trace formula Eq.\eqref{eq:traceformula}, in terms of the
signatures of the representations $\rho$ and $\rho\otimes\kan^{-2k}$.
Because $\dim\kan\!=\!1$, the latter may be determined by counting
the eigenvalue multiplicities of the matrices $\left(\textrm{-}1\right)^{k}\rho\!\sm 0{\textrm{-}1}10$
and $\exp\!\left(\textrm{-}\frac{2\pi\mathsf{i}}{3}\right)^{k}\rho\!\sm 0{\textrm{-}1}1{\textrm{-}1}$,
which is pretty straightforward, leading to the result summarized
in Table \ref{tab:sigs}.

\begin{table}[h]
\caption{\label{tab:sigs}The signature of $\rho\!\otimes\!\kan^{{\scriptscriptstyle -}2k}$ }
\begin{tabular}{|c|c|c|c|}
\hline 
$k\!\mod\,6$ & $\alpha\!\left(\rho\!\otimes\!\kan^{-2k}\right)$ & $\beta_{1}\!\left(\rho\!\otimes\!\kan^{-2k}\right)$ & $\beta_{2}\!\left(\rho\!\otimes\!\kan^{-2k}\right)$\tabularnewline
\hline
\hline 
0 & $\alpha$ & $\beta_{1}$ & $\beta_{2}$\tabularnewline
\hline 
1 & $d-\alpha$ & $\beta_{2}$ & $d-\beta_{1}-\beta_{2}$\tabularnewline
\hline 
2 & $\alpha$ & $d-\beta_{1}-\beta_{2}$ & $\beta_{1}$\tabularnewline
\hline 
3 & $d-\alpha$ & $\beta_{1}$ & $\beta_{2}$\tabularnewline
\hline 
4 & $\alpha$ & $\beta_{2}$ & $d-\beta_{1}-\beta_{2}$\tabularnewline
\hline 
5 & $d-\alpha$ & $d-\beta_{1}-\beta_{2}$ & $\beta_{1}$\tabularnewline
\hline
\end{tabular}

\end{table}

It follows that the sequence of $\gam k$-s is completely determined
by the signature of $\rho$, the first few values being summarized
in Table \ref{tab:gam}.

\begin{table}[h]
\caption{\label{tab:gam}The values of $\gam k$ for $0\leq k<6$}

\begin{tabular}{|c|c|}
\hline 
$k$ & $\gam k$\tabularnewline
\hline
\hline 
0 & $0$\tabularnewline
\hline 
1 & $\alpha+\beta_{1}+\beta_{2}-d$\tabularnewline
\hline 
2 & $\beta_{2}$\tabularnewline
\hline 
3 & $\alpha$\tabularnewline
\hline 
4 & $\beta_{1}+\beta_{2}$\tabularnewline
\hline 
5 & $\alpha+\beta_{2}$\tabularnewline
\hline
\end{tabular}
\end{table}

Inspection of Table \ref{tab:gam} reveals\emph{ }that\emph{ }$\gam k\!\geq\!0$
for $k\!>\!1$ (this can fail for $k\!=\!1$, a prime example being
the trivial representation $\rho_{0}$, for which $d\!=\!1$ and $\alpha\!=\!\beta_{1}\!=\!\beta_{2}\!=\!0$,
hence $\gam 1\!=\!-1$), and the relations $\gam 7\!=\!\gam 1\!+\! d\!=\!\gam 3\!+\!\gam 4$
and $\gam 5\!=\!\gam 2\!+\!\gam 3$. Thanks to Eq.\eqref{eq:gamdif},
these results generalize to $\gam{k+n}\!\geq\!\gam k$ for $n\!>\!1$,
and \begin{align}
\gam{k+7}+\gam k & =\gam{k+3}+\gam{k+4}\label{eq:gamsum2}\\
\gam{k+5}+\gam k & =\gam{k+3}+\gam{k+2}\,\,.\label{eq:gamsum1}\end{align}

The next important relation follows by considering the contragredient
representation $\dual{\rho}$, which assigns to each $\gamma\!\in\!\FD$
the transposed inverse of its representation operator $\rho\!\left(\gamma\right)$:
\begin{equation}
\dual{\rho}\!\left(\gamma\right)\!=\!\trans{\rho\!\left(\gamma^{-1}\right)}\,.\label{eq:contradef}\end{equation}
Since transposition does not change eigenvalues, while inversion inverts
them, it follows that the signature of $\dual{\rho}$ is given by
\begin{equation}
\dual d=d,\;\dual{\alpha}=\alpha,\;\dual{\beta_{1}}=\beta_{2},\;\dual{\beta_{2}}=\beta_{1}\,\,.\label{eq:dualsig}\end{equation}
Combining this with Table \ref{tab:gam} and the periodicity relation
Eq.\eqref{eq:gamperiod}, one arrives at the duality relation \begin{equation}
\gam k\!\left(\dual{\rho}\right)\!=\!\gam 1\!\left(\rho\right)-\gam{1-k}\!\left(\rho\right)\,.\label{eq:dual1}\end{equation}
On the other hand, by the definition of the contragredient \begin{equation}
\exp\!\left(2\pi\mathsf{i}\dual{\FB}\right)\!=\!\dual{\rho}\sm 1101\!=\!\trans{\rho}\!\sm 1{\textrm{-}1}0{\,1}\!=\!\exp\!\left(\textrm{-}2\pi\mathsf{i}\trans{\FB}\right)\,\,,\label{eq:duallambda}\end{equation}
if one denotes by $\dual{\FB}$ the exponent matrix of $\dual{\rho}$;
consequently, the sum $\trans{\FB}\!+\!\dual{\FB}$ has integer eigenvalues.
As a result, \begin{multline}
\tr{\ip{\trans{\FB}\!+\!\dual{\FB}\!-\!1}}\!=\!\tr{\left(\trans{\FB}\!+\!\dual{\FB}\!-\!1\right)}\!=\!\\
\tr{\,\left(\trans{\FB}\right)}+\tr{\,\dual{\FB}}-d\!=\! d\!-\!\alpha\!-\!\beta_{1}\!-\!\beta_{2}\!=\!-\gam 1\label{eq:dualtr}\end{multline}
according to Eqs.\eqref{eq:dualsig} and \eqref{eq:traceformula}.
But \begin{multline}
\kso\!\left(\dual{\rho}\right)\!=\!-\tr{\ip{1\!-\!\dual{\FB}}\!=\!-\tr{\ip{\trans{\FB}-\left(\trans{\FB}\!+\!\dual{\FB}\!-\!1\right)}}}\!=\\
\!-\!\tr{\ip{\trans{\FB}}}\!+\!\tr{\ip{\trans{\FB}\!+\!\dual{\FB}\!-\!1}}\!=\!-\tr{\ip{\FB}}\!-\!\gam 1\!\left(\rho\right)\,.\label{eq:dualkso}\end{multline}
proving the following supplement to Eq.\eqref{eq:dual1} \begin{equation}
\ks\!\left(\rho\right)+\kso\!\left(\dual{\rho}\right)\!=\!-\gam 1\!=\!\ks\!\left(\dual{\rho}\right)+\kso\!\left(\rho\right)\:.\label{eq:reciprocity}\end{equation}

The last major ingredient that we shall need is the observation that,
for a not necessarily irreducible even representation $\rho\!:\!\SL\!\rightarrow\!\gl{}V$
with finite image, a weight $0$ holomorphic form $\FA\!\in\!\mof 0{\rho}$
is a constant vector invariant under $\rho$; as a consequence, $\dim\mof 0{\rho}$
equals the multiplicity $\FC$ of the trivial representation in $\rho$,
and $\dim\cof 0{\rho}$ is always zero:\begin{equation}
\begin{aligned}\dim\mof 0{\rho} & =\FC\,,\\
\dim\cof 0{\rho}\, & =\,0\,.\end{aligned}
\label{eq:weight0}\end{equation}
To see this, note that any form $\FA\!\in\!\mof 0{\rho}$ is invariant
under the kernel $\ker\rho$ of $\rho$. This implies that, if $\mcr{\rho}$
denotes the surface uniformized by $\ker\rho$ and $\pi\!:\!\uhp\!\rightarrow\!\mcr{\rho}$
the associated natural projection, there exists a single valued map
$\hat{\FA}\!:\!\mcr{\rho}\!\rightarrow\! V$ such that $\FA\!=\!\hat{\FA}\circ\pi$.
Since $\FA$ is holomorphic on $\uhp$ and bounded at the cusp $\tau\!=\!\mathsf{i}\infty$,
the map $\hat{\FA}$ is holomorphic on all of $\mcr{\rho}$ (a finite
sheeted cover of the Riemann sphere $\mathbb{CP}^{1}$), and bounded
at all its cusps: by Liouville's theorem, it should be a constant
map. But this means that $\FA\!\in\!\mof 0{\rho}$ should be independent
of $\tau$, and this constant vector $\FA\!\in\! V$ has to satisfy
$\FA\!=\!\rho\!\left(\gamma\right)\FA$ for all $\gamma\!\in\!\FD$
because of Eq.\eqref{eq:modtrans}, hence it should be an invariant
vector of the representation $\rho$. Should $\FA$ be a cusp form,
it should vanish as $\tau\!\rightarrow\!\mathsf{i}\infty$, hence
it should vanish identically.

An immediate consequence of the above result is that, for $n\!>\!0$,
$\dim\mof{\textrm{-}2n}{\rho}\!\leq\!\dim\cof 0{\rho\!\otimes\!\kan^{\textrm{2}n}}\!=\!0$
by Eq.\eqref{eq:dimbound}, i.e. there are no holomorphic forms of
negative weight%
\footnote{Recall that we assume $\rho$ to have finite image.%
}. What is more, $\ks\!\leq\!0$ for irreducible and nontrivial $\rho$
(since such a $\rho$ has no invariant vectors), and $\kso\leq0$
for every $\rho$; combining the above with the reciprocity relation
Eq.\eqref{eq:reciprocity}, and noting that the contragredient $\dual{\rho}$
is irreducible and nontrivial whenever $\rho$ is, one gets the important
result that $\gam 1\!+\!\ks\!\geq\!0$ for all $\rho$, and \emph{$\gam 1\!+\!\kso\!\geq\!0$
}for irreducible and nontrivial $\rho$.

Let's now suppose that $\rho$ is even irreducible of dimension $d\!>\!1$.
In this case $\gam 1\!\geq\!-\kso\!\geq\!0\!=\!\gam 0$ by the above.
But if $\rho$ is irreducible of dimension $d\!>\!1$, then the same
is true of its tensor product with any representation of dimension
$1$, hence $\gam 1\!\left(\rho\!\otimes\!\kan^{\textrm{-2}k}\right)\!\geq\!0$
for all $k$. But $\gam 1\!\left(\rho\!\otimes\!\kan^{\textrm{-2}k}\right)\!=\!\gam{k+1}-\gam k$
according to Eq.\eqref{eq:gamdif}, leading to the conclusion that
in this case the $\gam k$-s form an increasing sequence: \begin{equation}
\ldots\leq\gam{-1}\leq\gam 0\!=\!0\leq\gam 1\leq\gam 2\leq\ldots\label{eq:ineq}\end{equation}
(note that this fails for $d\!=\!1$). Combining this result with
\emph{$\gam 1\!+\!\kso\!\geq\!0$}, one gets that for an irreducible\emph{
$\rho$ }with $d\!>\!1$ the inequality $\kso\!+\!\gam k\!\geq\!0$
holds for all $k\!>\!0$; a simple case by case check shows that it
does also hold for all nontrivial $\rho$ with $d\!=\!1$.

Finally, putting everything together, and taking into account Eq.\eqref{eq:dimforred},
we get that for an even, not necessarily irreducible representation
$\rho$\begin{equation}
\dim\mof{2k}{\rho}=\begin{cases}
0 & \mathrm{if}\; k<0;\\
\FC & \mathrm{if}\; k=0;\\
\ks+\gam k & \mathrm{if}\; k>0,\end{cases}\label{eq:dimformula2}\end{equation}
 and \begin{equation}
\dim\cof{2k}{\rho}=\begin{cases}
0 & \mathrm{if}\; k\leq0;\\
\kso+\gam 1+\FC & \mathrm{if}\; k=1;\\
\kso+\gam k & \mathrm{if}\; k>1.\end{cases}\label{eq:dimformula3}\end{equation}

It follows from Eqs.\eqref{eq:dimformula2} and \eqref{eq:dimformula3},
combined with the duality relations Eqs.\eqref{eq:dual1} and \eqref{eq:reciprocity},
that \begin{equation}
\dim\mof{2k}{\rho}+\dim\cof{12n+2-2k}{\dual{\rho}}\!=\! nd\label{eq:dual2}\end{equation}
for positive integers $k$ and $n$ such that $k\!<\!6n$, expressing
the duality between holomorphic and cusp forms, and \begin{equation}
\begin{aligned}\ks\!\left(\rho\right) & =\dim\mof 0{\rho}-\dim\cof 2{\dual{\rho}}\\
\kso\!\left(\rho\right) & =\dim\cof 0{\rho}-\dim\mof 2{\dual{\rho}}\,.\end{aligned}
\label{eq:lambdacohom}\end{equation}

Let's take a look at the classical case, when $\rho$ is the identity
representation $\rho_{0}$. In this case $\ks\!=\!\FC\!=\!1$, $\kso\!=\!\alpha\!=\!\beta_{1}\!=\!\beta_{2}\!=\!0$,
and for $k\!<\!6$ all $\gam k$-s are zero except for $\gam 1=-1$.
This leads to the well-known result \begin{equation}
\dim\mathsf{M}_{2k}=\begin{cases}
\left[\frac{k}{6}\right] & \mathrm{if}\:\Mod k16\:,\\
\left[\frac{k}{6}\right]+1 & \mathrm{otherwise},\end{cases}\label{eq:classicaldim}\end{equation}
and $\dim\mathsf{S}_{2k}\!=\!\max\!\left(0,\dim\mathsf{M}_{2k}\!-\!1\right)$.
As we can see, the classical case is somewhat exceptional because
of the existence of an invariant vector for $\rho_{0}$.

To finish, let's consider the weight distribution of the generators
of $\mof{}{\rho}\!=\!\oplus_{k}\mof k{\rho}$ and $\cof{}{\rho}\!=\!\oplus_{k=0}^{\infty}\cof k{\rho}$,
considered as free modules over the ring $\mathbb{C}\!\left[E_{4},E_{6}\right]$
of holomorphic scalar modular forms, cf. Eq.\eqref{eq:modring}. The
first step is to compute the Hilbert-Poincaré series \begin{multline}
\Gmof{\rho}z\!=\!\sum_{k=0}^{\infty}\dim\mof{2k}{\rho}z^{2k}\!=\!\FC+\sum_{k=1}^{\infty}\left(\ks\!+\!\gam k\right)z^{2k}=\\
=\FC+\frac{\ks z^{2}}{1-z^{2}}+\frac{\gam 1z^{2}\!+\!\gam 2z^{4}\!+\!\left(\gam 3\!-\!\gam 1\right)z^{6}\!+\!\left(\gam 4\!-\!\gam 2\!-\!\gam 1\right)z^{8}}{\left(1-z^{4}\right)\left(1-z^{6}\right)}\,\label{eq:genfun}\end{multline}
and \begin{multline}
\Gcof{\rho}z\!=\!\sum_{k=0}^{\infty}\dim\cof{2k}{\rho}z^{2k}\!=\!\left(\gam 1\!+\!\kso\!+\!\FC\right)\! z^{2}\!+\!\sum_{k=2}^{\infty}\!\left(\kso\!+\!\gam k\right)\! z^{2k}\\
=\FC z^{2}+\frac{\kso z^{2}}{1-z^{2}}+\frac{\gam 1z^{2}\!+\!\gam 2z^{4}\!+\!\left(\gam 3\!-\!\gam 1\right)z^{6}\!+\!\left(\gam 4\!-\!\gam 2\!-\!\gam 1\right)z^{8}}{\left(1-z^{4}\right)\left(1-z^{6}\right)}\,,\label{eq:genfunB}\end{multline}
where we have used the relations Eqs.\eqref{eq:gamsum2} and \eqref{eq:gamsum1}
to sum up the power series $\sum_{k=0}^{\infty}\gam kz^{2k}$. From
Eqs.\eqref{eq:genfun} and \eqref{eq:genfunB} one reads off the weight
distribution of the generators as the coefficients of the polynomials
$\left(1\!-\! z^{4}\right)\!\left(1\!-\! z^{6}\right)\Gmof{\rho}z$
and $\left(1\!-\! z^{4}\right)\!\left(1\!-\! z^{6}\right)\Gcof{\rho}z$,
the results being tabulated in Table \ref{Flo:genweights}. 

\begin{table}[h]
\caption{Weight distribution of free generators for even $\rho$.}
\label{Flo:genweights}\begin{tabular}{|c|l|r|}
\hline 
weight & $\mof{}{\rho}$  &  $\cof{}{\rho}$\tabularnewline
\hline
\hline 
0 & $\FC$ & $0$\tabularnewline
\hline 
2 & $\gam 1\!+\!\ks$ & $\gam 1\!+\!\kso\!+\!\FC$\tabularnewline
\hline 
4 & $\gam 2\!+\!\ks\!-\!\FC$ & $\gam 2\!+\!\kso$\tabularnewline
\hline 
6 & $\gam 3\!-\!\gam 1\!-\!\FC$ & $\gam 3\!-\!\gam 1\!-\!\FC$\tabularnewline
\hline 
8 & $\gam 6\!-\!\gam 5\!-\!\ks$ & $\gam 6\!-\!\gam 5\!-\!\kso\!-\!\FC$\tabularnewline
\hline 
10 & $\FC\!-\!\ks$ & $-\kso$\tabularnewline
\hline 
12 & $0$ & $\FC$\tabularnewline
\hline
\end{tabular}
\end{table}

Note that in both cases we have a total of $\gam 6\!=\! d$ generators,
in accord with the fact that these are free modules of rank $d$.
Finally, the underlying duality Eq.\eqref{eq:dual2} between cusps
forms for $\rho$ and holomorphic forms for its contragredient is
elegantly expressed by the relation \begin{equation}
\Gcof{\dual{\rho}}z\!=\! z^{2}\Gmof{\rho}{z^{{\scriptscriptstyle -}1}}\,\label{eq:HPdual}\end{equation}
between the respective Hilbert-Poincaré series%
\footnote{This has to be interpreted as a relation between the corresponding
rational expressions to which these power series sum up.%
}.

\section{Odd representations\label{sec:odd}}

\global\long\def\dks{\dot{\boldsymbol{\uplambda}}_{{\scriptscriptstyle +}}}
\global\long\def\dkso{\dot{\boldsymbol{\uplambda}}_{{\scriptscriptstyle -}}}
\global\long\def\dgam#1{\dot{\upgamma}_{#1}}

Let's consider an odd representation $\rho\!:\!\SL\!\rightarrow\!\gl{}V$
with finite image. In this case, nontrivial forms have odd weight,
and the representation $\dot{\rho}\!=\!\rho\otimes\kan^{-1}$ is even,
and still has finite image: we shall denote its exponent matrix by
$\dot{\FB}$, and introduce the notations \begin{equation}
\dgam k\!\left(\rho\right)\!=\!\gam k\!\left(\dot{\rho}\right)\label{eq:dgamdef}\end{equation}
and\begin{equation}
\begin{aligned}\dks\!\left(\rho\right) & =\!\tr{\!\ip{\dot{\FB}+\frac{1}{12}}}\\
\dkso\!\left(\rho\right) & =\!-\!\tr{\!\ip{\frac{11}{12}\!-\!\dot{\FB}}}\,.\end{aligned}
\label{eq:lamdefodd}\end{equation}
Note that comments similar to those after Eq.\eqref{eq:dimfor4} apply
in this case as well, e.g. $\dgam{k+6}\!=\!\dgam k\!+\! d$. Moreover,
since $\ip{y\!-\! x}\!+\!\ip x\!\leq\!0$ for all $x\!\in\!\mathbb{R}$
and $y\!<\!1$, one has \begin{equation}
\ks\!\left(\dot{\rho}\right)\!\leq\!\dkso\!\left(\rho\right)\!\leq\!\dks\!\left(\rho\right)\,,\label{eq:lambdaineq}\end{equation}
the second inequality following from the observation that $\dks\!-\!\dkso$
equals the number of invariant vectors of $\rho\!\sm 1101$, hence
it can't be negative. Recalling from the discussion following Eq.\eqref{eq:ineq}
that $\gam k\!+\!\ks\!\geq\!0$ for even representations and positive
$k$, we get from Eq.\eqref{eq:lambdaineq} the following inequality
for odd $\rho$ and $k\!>\!0$: \begin{equation}
\dks\!\left(\rho\right)\!+\!\dgam k\!\left(\rho\right)\!\geq\!\dkso\!\left(\rho\right)\!+\!\dgam k\!\left(\rho\right)\!\geq\!\ks\!\left(\dot{\rho}\right)\!+\!\gam k\!\left(\dot{\rho}\right)\!\geq\!0\,.\label{eq:posdim}\end{equation}

There is, however, an important difference with respect to the case
of even $\rho$: namely, since $\!\dot{\left(\dual{\rho}\right)}\!=\!\dual{\left(\dot{\rho}\right)}\otimes\kan^{-2}$
as a consequence of $\dual{\kan}\!=\!\kan^{-1}$, for odd representations
the duality relations take the form\begin{equation}
\dgam k\!\left(\dual{\rho}\right)\!=\!-\dgam{-k}\!\left(\rho\right)\,\label{eq:dgamrecip}\end{equation}
and \begin{equation}
\dks\!\left(\dual{\rho}\right)=-\dkso\!\left(\rho\right)\,,\label{eq:reciprocityodd}\end{equation}
to be compared with Eqs.\eqref{eq:dual1} and \eqref{eq:reciprocity}. 

The next question concerns the minimal possible weight of a holomorphic
or cusp form for odd $\rho$. The case of $\rho\!=\!\kan$ shows that
there could exist weight 1 cusp forms, the prime example being the
form $\un\!\in\!\cof 1{\kan}$ entering the weight-shifting map from
Section \ref{sec:general}. On the other hand, there can be no nontrivial
holomorphic forms of negative weights, for Eq.\eqref{eq:dimbound}
and Eq.\eqref{eq:weight0} imply that for $k\!\geq\!0$ \begin{equation}
\dim\cof{-\left(2k+1\right)}{\rho}\leq\dim\mof{-\left(2k+1\right)}{\rho}\leq\dim\cof 0{\rho\otimes\kan^{2k+1}}\!=\!0\,.\label{eq:lowweightbound}\end{equation}

For irreducible $\rho$, an argument paralleling the one leading to
Eq.\eqref{eq:dimfor4} gives at once\begin{equation}
\begin{aligned}\dim\mof{2k+1}{\rho} & =\max\!\left(0,\dks\!+\!\dgam k\,\right)\,,\\
\dim\cof{2k+1}{\rho} & =\max\!\left(0,\dkso\!+\!\dgam k\,\right)\,.\end{aligned}
\label{eq:dimfor5}\end{equation}
Taking into account the above discussed inequalities, one concludes
that for odd irreducible $\rho$\begin{equation}
\dim\mof{2k+1}{\rho}=\begin{cases}
0 & \mathrm{if}\; k<0;\\
\dks+\dgam k & \mathrm{if}\; k>0,\end{cases}\label{eq:dimformula2odd}\end{equation}
 and \begin{equation}
\dim\cof{2k+1}{\rho}=\begin{cases}
0 & \mathrm{if}\; k<0;\\
\dkso+\dgam k & \mathrm{if}\; k>0,\end{cases}\label{eq:dimformula3odd}\end{equation}
which should be supplemented with the relations $\dim\mof 1{\rho}\!=\!\max\!\left(0,\dks\right)$
and $\dim\cof 1{\rho}\!=\!\max\!\left(0,\dkso\right)$. Combining
the above results with Eqs. \eqref{eq:dgamrecip} and \eqref{eq:reciprocityodd},
one obtains the duality relation \begin{equation}
\dim\mof{2k+1}{\rho}+\dim\cof{12n+1-2k}{\dual{\rho}}\!=\! nd\,,\label{eq:dualodd}\end{equation}
valid for positive integers $k$ and $n$ such that $k\!<\!6n$, and\begin{equation}
\begin{aligned}\dks= & \dim\mof 1{\rho}-\dim\cof 1{\dual{\rho}}\\
\dkso= & \dim\cof 1{\rho}-\dim\mof 1{\dual{\rho}}\,,\end{aligned}
\label{eq:dlambdacohom}\end{equation}
explaining the meaning of the parameters $\dks$ and $\dkso$. Note
that, taking into account Eq.\eqref{eq:dimforred}, this last result
holds for arbitrary, not necessarily irreducible odd representations,
because both sides of the equalities hold for each irreducible constituent
separately and the relevant quantities are additive. What is more,
by an analogue argument the same is true of Eqs.\eqref{eq:dimformula2odd}
and \eqref{eq:dimformula3odd}, which should be regarded as the final
form of the dimension formula for odd representations%
\footnote{Note that this allows us to dispense with irreducible decomposition
of $\rho$, since the dimensions can be expressed solely in terms
of some global parameters characterizing $\rho$. This is to be contrasted
with the weight 1 case, where one needs the explicit knowledge of
the irreducible decomposition to be able to compute the relevant dimensions.%
}.

As an example, let's consider the odd representation $\kan$. Since
$\dot{\kan}\!=\!\rho_{0}$, one has $\dks\!\left(\kan\right)\!=\!\dkso\!\left(\kan\right)\!=\!1$,
and $\dgam k\!=\!0$ for $0\leq k\!<\!7$, except for $\dgam 1\!=\!-1$.
It follows that $\dim\mof 1{\kan}\!=\!\dim\cof 1{\kan}\!=\!1$ (and
indeed, $\un\!\in\!\cof 1{\kan}$ as remarked before), while $\dim\mof 3{\kan}\!=\!0$.
From Eq.\eqref{eq:dlambdacohom}, we see that $\dim\mof 1{\dual{\kan}}\!=\!0$.

Finally, let's discuss the weight distribution of the generators of
$\mof{}{\rho}\!=\!\oplus\mof k{\rho}$ and $\cof{}{\rho}\!=\!\oplus\cof k{\rho}$,
considered as modules over the ring $\mathsf{M}\!=\!\mathbb{C}\!\left[E_{4},E_{6}\right]$
of scalar holomorphic forms. An argument completely parallel to that
leading to Eq.\eqref{eq:genfun}, but based on the results relevant
to odd representations, gives Table \ref{Flo:genweightsodd} for the
weight distribution of the generators (note that we again have a total
of $\dgam 6\!=\! d$ independent generators). The duality relation
Eq.\eqref{eq:HPdual} goes over verbatim to the odd case.

\begin{table}[h]
\caption{Weight distribution of free generators for odd $\rho$.}
\label{Flo:genweightsodd}\begin{tabular}{|c|l|r|}
\hline 
weight & $\mof{}{\rho}$  &  $\cof{}{\rho}$\tabularnewline
\hline
\hline 
1 & $\dim\mof 1{\rho}$ & $\dim\cof 1{\rho}$\tabularnewline
\hline 
3 & $\dgam 1\!+\!\dks$ & $\dgam 1\!+\!\dkso$\tabularnewline
\hline 
5 & $\dgam 2\!+\!\dks\!-\!\dim\mof 1{\rho}$ & $\dgam 2\!+\!\dkso\!-\!\dim\cof 1{\rho}$\tabularnewline
\hline 
7 & $\dgam 3\!-\!\dgam 1\!-\!\dim\mof 1{\rho}$ & $\dgam 3\!-\!\dgam 1\!-\!\dim\cof 1{\rho}$\tabularnewline
\hline 
9 & $\dgam 6\!-\!\dgam 5\!-\!\dks$ & $\dgam 6\!-\!\dgam 5\!-\!\dkso$\tabularnewline
\hline 
11 & $\dim\mof 1{\rho}\!-\!\dks$ & $\dim\cof 1{\rho}\!-\!\dkso$\tabularnewline
\hline
\end{tabular}
\end{table}

\section{Outlook}

\noindent We have investigated spaces of vector-valued holomorphic
and cusp forms of integer weight for finite dimensional representations
of the modular group $\FD\!=\!\SL$ having finite image, and have
obtained explicit expressions, see Eqs.\eqref{eq:dimformula2}, \eqref{eq:dimformula3},
\eqref{eq:dimformula2odd} and \eqref{eq:dimformula3odd}, for the
dimension of these spaces. Based on these results, we have described
the weight distribution of the generators of the module of holomorphic
and cusp forms, and the duality, most elegantly expressed by Eq.\eqref{eq:HPdual},
relating cusp forms with holomorphic forms for the contragredient.

It goes without saying that the results presented here agree completely
with those of \cite{Skoruppa1984}: to see this, one has to rewrite
the quantities appearing in \cite{Skoruppa1984} in terms of those
of the present paper, which is a straightforward job using Eq.\eqref{eq:sigtrace}.
The advantage of our approach is that it doesn't only give us the
actual dimensions, but it does also provide an effective procedure
for computing explicit bases, by solving the relevant system of linear
relations, as described in the argument leading to Eq.\eqref{eq:dimfor1}.
As an extra bonus, we get a much better control over the quantities
involved, making it easier to recognize relations like Eq.\eqref{eq:HPdual}.

Several possible generalizations offer themselves at once. First,
one could try to generalize the theory from integer weight to half-integer,
or even arbitrary real weights: this would necessitate the consideration
of suitable projective representations of $\SL$, making the whole
story a bit more complicated than in the integer weight case. Next,
one could contemplate the possibility to dispense of the finite image
requirement: this could result in severe difficulties, both analytic
(the surface uniformized by the kernel would not be anymore a finite
sheeted cover of the modular curve, allowing for holomorphic forms
of negative weight) and algebraic (the possibility of reducible but
indecomposable representations could lead to non-semisimple exponent
matrices), but is certainly a most interesting issue to be dealt with,
since this is the case relevant for logarithmic conformal theories
(to be contrasted with the finite image case relevant for rational
theories). Finally, an obvious generalization, making direct contact
with classical knowledge, would be to consider forms for (finite index)
subgroups of $\SL$. To sum up, there are many interesting questions
left open for future investigations.

\bibliographystyle{plain}


\end{document}